\documentclass[11pt,a4paper]{amsart}

\usepackage{amsmath}
\usepackage{amssymb}
\usepackage{enumerate}
\usepackage{hyperref}

\newtheorem{theorem}{Theorem}[section]
\newtheorem{proposition}[theorem]{Proposition}
\newtheorem{lemma}[theorem]{Lemma}

\newtheorem{remark}[theorem]{Remark}

\newtheoremstyle{questionstyle}
  {}{}                  
  {\itshape}            
  {}                   
  {\bfseries}           
  {.}                   
  { }                   
  {\ifblank{#3}{\thmname{#1}\thmnumber{ #2}}{#3}}

\theoremstyle{questionstyle}
\newtheorem{question}{question}

\newcommand{\f}{\mathbf{f}}
\newcommand{\e}{\mathbf{e}}

\newcommand{\B}{\mathcal{B}}
\newcommand{\A}{\mathcal{A}}

\newcommand{\N}{\mathbb{N}}

\newcommand{\R}{\mathbb{R}}
\newcommand{\F}{\mathbb{F}}
\newcommand{\K}{\mathbb{K}}
\renewcommand{\L}{\mathbb{L}}
\newcommand{\E}{\mathbb{E}}

\newcommand{\m}{\mathfrak{m}}
\newcommand{\p}{\mathfrak{p}}

\newcommand{\M}{\mathrm{M}}
\newcommand{\End}{\mathrm{End}}
\newcommand{\colen}{\ell}
\newcommand{\vol}{\mathrm{Vol}}
\newcommand{\diff}{\mathrm{d}}
\newcommand{\Ann}{\mathrm{Ann}}
\renewcommand{\ker}{\mathrm{Ker}}
\newcommand{\ini}{\mathrm{in}}
\newcommand{\LC}{\mathrm{LC}}
\newcommand{\LM}{\mathrm{LM}}

\newcommand{\Span}{\mathrm{Span}}

\begin{document}

    \title[Algebras generated by commuting matrices]{Linear bounds and sharp exponential growth for commuting matrix algebra}

    \author{Jiamei Li}
    \address{SCHOOL OF MATHEMATICAL SCIENCES, HEBEI NORMAL UNIVERSITY, SHIJIAZHUANG, HEBEI 050024, CHINA.}
    \email{lijm\_03@163.com}

    \author{Caihua Luo}
    \address{SCHOOL OF SCIENCE AND ENGINEERING, THE CHINESE UNIVERSITY OF HONG KONG, SHENZHEN,
GUANGDONG 518172, CHINA.}
    \email{luocaihua@cuhk.edu.cn}

    \author{Xu Song}
    \address{SCHOOL OF MATHEMATICAL SCIENCES, HEBEI NORMAL UNIVERSITY, SHIJIAZHUANG, HEBEI 050024, CHINA.}
    \email{songx@hebtu.edu.cn}

    \author{Yilin Zeng}
    \address{SCHOOL OF SCIENCE AND ENGINEERING, THE CHINESE UNIVERSITY OF HONG KONG, SHENZHEN,
GUANGDONG 518172, CHINA.}
    \email{223015132@link.cuhk.edu.cn}

    %\date{}


    \begin{abstract}
        For any $m$ pairwise commuting matrices $\{A_i\}_{i=1}^{m}\subset \M_n(\F)$, we prove $\dim_\F \F[A_1,\cdots, A_m]\leq (2^m-1)(n-1)+1$, answering Bergman's question on the existence of a linear bound for a fixed number $m$.
    \end{abstract}

    \maketitle

    \section{Introduction}

        Let $\F$ be a field. For pairwise commuting matrices $A_1,\cdots, A_m\in\M_n(\F)$, this paper establishes linear upper bounds for the dimension of $\F[A_1,\cdots, A_m]$.

        \subsection*{Gerstenhaber problem}
        For a single matrix $A_1$, the Cayley--Hamilton theorem says that the dimension of $\F[A_1]$ is at most $n$. For two commuting matrices $A_1, A_2$, the dimension of $\F[A_1, A_2]$ is again at most $n$, which was proved by Motzkin and Taussky \cite{MR86781} and, independently, by Gerstenhaber \cite{MR132079}, both arguments resting on the irreducibility of the variety of commuting pairs. Linear-algebraic and commutative-algebraic proofs were found later \cite{MR1064891, MR1088666, MR1064892, bergman2013commutingmatricesmodulesartinian}. For four or more matrices, the corresponding bound fails \cite{MR132079}. Correspondingly, the variety of commuting $m$-tuples of $n\times n$ matrices is reducible when $m,n\geq 4$ \cite{MR132079, MR1199042}.

        The case of three matrices lies between these two regimes and remains open.

        \begin{question}[Gerstenhaber Problem]
            Is $\dim_\F \F[A,B,C]\leq n$ for any pairwise commuting $A,B,C\in\M_n(\F)$? 
        \end{question}

        We refer to \cite{MR3013206, MR3278252} for surveys. For $n\leq 10$, the answer is affirmative: in this range, the variety of commuting triples of $n\times n$ matrices over an algebraically closed field is irreducible \cite{MR1753173, MR1823346, MR2073906, MR2181883, MR2401636, MR2921710}, and Gerstenhaber's argument for pairs then applies. For $n\geq 29$, the variety of commuting triples is reducible \cite{MR1823346, MR3201099} (see also \cite{MR4445543} for its components), so the geometric approach cannot settle this problem for all $n$. Affirmative answers are known for various special classes of triples, for example, when one of the matrices has square zero \cite{MR1823346} or small nullity \cite{MR1753173}, and for several combinatorially defined classes \cite{MR3863478, MR4143269, MR4826769,cherny2025combinatorialalgebrasgeneratedcommuting}. Holbrook and O'Meara \cite{holbrook2020computingstrategyprogramsresolve} have proposed a computational strategy for resolving the problem. 

        \subsection*{Bergman question}
         It is well-known that every commutative subalgebra of $\M_n(\F)$ has dimension at most $\lfloor \frac{n^2}{4}\rfloor+1$, by theorems of Schur \cite{MR1580676} and Jacobson \cite{MR10540} (see also \cite{MR1615548}), and this bound is sharp; the extremal algebras, however, require on the order of $n^2$ generators. For a fixed number of generators, subquadratic bounds arose in the representation theory of $p$-adic groups, in connection with Bernstein's uniform admissibility theorem \cite{Ber74}. It states that $m$ pairwise commuting complex $n\times n$ matrices generate an algebra of dimension at most $n^{2-2^{1-m}}$ (see \cite[\S 4.10]{MR425030}\cite[Proposition 20]{Bern92}). 

        In the other direction, Bergman \cite[\S 6]{bergman2013commutingmatricesmodulesartinian} considered the supremum $r_m$ of the ratios $\dim_\F A/\dim_\F V$ over all commutative subalgebras $A\subset\End_\F(V)$ generated by $m$ elements. Explicitly, define
        \begin{equation*}
            c_m(\F)=
            \sup_{\substack{n\ge1,\ A_1,\cdots,A_m\in\M_n(\F)\\
            A_iA_j=A_jA_i}}
            \frac{\dim_\F \F[A_1,\cdots,A_m]}{n},
            \qquad r_m=\sup_\F c_m(\F).
        \end{equation*}
        He asked the following.
        \begin{question}[Bergman Question]
            Let $m\geq 3$. Is this $r_m$ finite?
        \end{question}

        All the bounds above grow faster than linearly in $n$, so they do not address Bergman's question.

        \subsection*{Linear Bound}
        Our main result answers the Bergman question affirmatively. 
        
        \begin{theorem}\label{thm: linear}
        Let $\F$ be any field and let $m,n\ge1$. For any $A_1,\ldots,A_m\in\M_n(\F)$ commute pairwise, then
        \begin{enumerate}[(i)]
        	\item $\dim_\F \F[A_1,\ldots,A_m]\le (2^m-1)(n-1)+1$.
        	\item $\lim\limits_{m\to\infty}c_m(\F)^{\frac{1}{m}}=
        	\lim\limits_{m\to\infty}r_m^{\frac{1}{m}}=2$.
        	\item if $\mathrm{char}(\F)=0$ (resp. a prime $p$), then $c_m(\F)=c_m(\mathbb{Q})$ (resp. $c_m(\F_p)$). 	
        \end{enumerate}
        \end{theorem}

        \begin{remark}
        In fact, while obtaining Theorem \ref{thm: linear} in this paper, we also get a sharper bound for $r_m$ for every $m$ as follows.
        \begin{equation*}
            r_m\leq \frac{\Gamma_{\!E}(m)}{\Gamma_{\!E}\left(\frac{m+1}{2}\right)^2}.
        \end{equation*}
        In particular, for the Gerstenhaber problem, we have $\dim_\mathbb{F} \mathbb{F}[A_1,A_2,A_3]\leq 2n-1$, i.e.,  $r_3\leq 2$. The detailed proof of $r_m$ will appear in a subsequent work as it contains new ideas.
    \end{remark}

\emph{\bf \underline{Sketch of Proof of (i) and (ii)}}

(i) For the upper bound, the first ingredient is the maximal-rank argument. Put $V=\F^n$ and $R=\F[A_1,\ldots, A_m]$, and choose $v$ so that $W=Rv$ has maximal dimension, we have $\Ann_R(v)$ sends $V$ into $W$ and annihilates $W$. Write $P=\F[x_1,\ldots,x_m]$, and let $Q,H,J$ be the annihilators in $P$ of $V,W,V/W$, respectively. Then
\[
H\subseteq J,\qquad HJ\subseteq Q,\qquad
\dim_\F(P/H)=\dim_\F W.
\]
The second ingredient is the colength inequality
\[
\bigl(\dim_\F(P/KL)\bigr)^{1/m}
 \le \bigl(\dim_\F(P/K)\bigr)^{1/m}
       +\bigl(\dim_\F(P/L)\bigr)^{1/m}.
\]
If $r=\dim_\F W$ and $s=\dim_\F(P/J)$, these ingredients yield
\[
\dim_\F R\le(r^{1/m}+s^{1/m})^m\le(2^m-1)r+s,
\]
because $s\le r$. Induction applies to the algebra induced on $V/W$.

(ii) The exponential growth rate follows from constructing a lower bound. Write $m=p+q$ with $p,q\ge2$. From the coefficients of the product of two polynomials of degrees $p-1$ and $q-1$, together with one additional bilinear monomial, we obtain a graded algebra $\A$ with $m$ degree-one generators. Its ideal of relations is generated by a single homogeneous polynomial of total degree
\[
e=\binom{m-2}{p-1}.
\]
Formal shift operators then give a faithful action of $\A/\A_{>d}$ on a space of dimension $\binom{d+m-1}{m-1}$. The ratio of the dimensions of the algebra to this space tends to $e$ as $d\to\infty$. Balancing $p-1$ and $q-1$ shows that $c_m(\F)\geq \binom{m-2}{\lfloor\frac{m-2}{2}\rfloor}$, and the central binomial coefficient has exponential rate $2$.

    \section{Proof of Theorem \ref{thm: linear}(i)}

    \subsection{Colength inequality ingredient}

        Let $P=\F[x_1,\cdots,x_m]$. For an ideal $I$ in $P$, we define
        \begin{equation*}
            \colen(I) := \dim_\F(P/I).
        \end{equation*}
        We say $I$ has finite colength if $\colen(I)$ is finite. 
        If $I$ is a monomial ideal, i.e., $I$ is generated by some monomials $x^\alpha=x_1^{\alpha_1}\cdots x_m^{\alpha_m}$ with an $m$-tuple $\alpha=(\alpha_1,\cdots,\alpha_m)\in\N^m$, we define
        \begin{equation*}
            \Gamma(I) = \bigcup_{x^\alpha\in I}(\alpha+\R_{\geq 0}^m).
        \end{equation*}

        \begin{lemma}\label{lem: colength}
            An ideal $I$ of $P$ has finite colength if and only if there exist $a_i\in\N_+$ with $x_i^{a_i}\in I$ for each $i=1,\cdots,m$.
        \end{lemma}
        \begin{proof}
            If $x_i^{a_i}\not\in I$ for all $a_i\geq 1$. Then the images of $1,x_i,x_i^2,\cdots$ in $P/I$ are linearly independent, so $\colen(I)=\infty$. Conversely, assume $x_i^{a_i}\in I$ for each $i$, and let $x_i^\alpha\not\in I$. Then $\alpha_i\leq a_i-1$ for each $i$ since $x_i^{a_i}\mid x^\alpha$ implies $x^\alpha\in I$. So the monomials outside $I$ are contained in $\{x_1^{\alpha_1}\cdots x_m^{\alpha_m}:0\leq \alpha_i\leq a_i-1\}$, a set of at most $\prod\limits_{i} a_i$ elements. Therefore, $\colen(I)\leq \prod\limits_i a_i<\infty$.
        \end{proof}
        
        \begin{lemma}\label{lem:inequality}
            Let $K, L$ be two monomial ideals of $P$ that have finite colength. Then
            \begin{enumerate}[(i)] 
                \item $\Gamma(K)$ is closed.
                \item $\colen(K) = \vol(\R_{\geq 0}^m\setminus\Gamma(K))$.
                \item $KL$ has finite colength.
                \item $\Gamma(KL) = \Gamma(K) + \Gamma(L)$.
            \end{enumerate}
        \end{lemma}
        \begin{proof}
            (i) By Lemma \ref{lem: colength}, there exist some $a_i\in\N_+$ such that $x_i^{a_i}\in K$. Let $e_j=(\delta_{1j},\cdots,\delta_{mj})$ and $\delta_{ij}$ is the Kronecker delta which is $1$ if $i=j$, $0$ otherwise. Then
            \begin{equation*}
                \Gamma(K) = \left(\bigcup_{j=1}^m(a_je_j+\R_{\geq 0}^m)\right) \cup \left(\bigcup_{\substack{x^\alpha\in K\\ \alpha_i\leq a_i-1\;\text{for all i}}}(\alpha+\R_{\geq 0}^m)\right)
            \end{equation*}
            since if $\alpha_j\geq a_j$, then $\alpha\geq a_je_j$ coordinate-wise, and hence $\alpha+\R_{\geq0}^m\subset a_je_j+\R_{\geq 0}^m$, otherwise $\alpha\in\prod\limits_{i}\{0,1,\cdots,a_i-1\}$. Then $\Gamma(K)$ is a finite union of closed subsets. Therefore, $\Gamma(K)$ is closed.

            (ii) For $\alpha\in\R_{\geq 0}^m$, we claim $\alpha\in\Gamma(K)$ if and only if $x^{\lfloor\alpha\rfloor}\in K$ where $\lfloor\alpha\rfloor = (\lfloor\alpha_1\rfloor,\cdots,\lfloor\alpha_m\rfloor)$. If $\alpha\in\Gamma(K)$, then $\alpha =\beta+u$ with $x^\beta\in K$ and $u\in\R_{\geq 0}^m$, so $\beta_i\leq \alpha_i$ for each $i$. Since $\beta_i$ are integers, $\beta_i\leq \lfloor\alpha_i\rfloor$, and $x^{\lfloor\alpha\rfloor} = x^{\beta}x^{\lfloor\alpha\rfloor-\beta}\in K$. Conversely, if $x^{\lfloor\alpha\rfloor}\in K$, then $\alpha\in \lfloor\alpha\rfloor+\R_{\geq0}^m\subset\Gamma(K)$. 
            
            It follows that 
            \begin{equation*}
                \R_{\geq0}^m\setminus\Gamma(K) = \bigsqcup_{x^\alpha\notin K} (\alpha +[0,1)^m),
            \end{equation*}
            since for any $\alpha\in \R_{\geq 0}^m$, $\alpha\in\lfloor\alpha\rfloor+[0,1)^m$, this $m$-cube lies in the complement if and only if $x^{\lfloor\alpha\rfloor}\notin K$, and distinct $\lfloor\alpha\rfloor+[0,1)^m$ are disjoint.
            By Lemma \ref{lem: colength}, there exist finitely many such $m$-cubes, each of volume $1$, so the right-hand side is measurable and 
            \begin{equation*}
                \vol(\R_{\geq 0}^m\setminus\Gamma(K)) = \#(\{x^\alpha:x^{\alpha}\notin K\}) = \colen(K).
            \end{equation*}

            (iii) By Lemma \ref{lem: colength}, there exist $a_i, b_i\in \N_+$ such that $x_i^{a_i}\in K$ and $x_i^{b_i}\in L$. Then $x_i^{a_i+b_i} = x_i^{a_i}x_i^{b_i}\in KL$ for each $i$. By Lemma \ref{lem: colength}, $KL$ has finite colength.

            (iv) If $x^\alpha\in K $ and $x^\beta\in L$, then $x^{\alpha+\beta} = x^{\alpha}x^{\beta}\in KL$, so
            \begin{equation*}
                (\alpha+\R_{\geq0}^m) + (\beta+\R_{\geq 0}^m) = (\alpha+\beta) + \R_{\geq 0}^m \subset \Gamma(KL).
            \end{equation*}
            It follows that $\Gamma(K)+\Gamma(L)\subset \Gamma(KL)$. Conversely, let $\gamma\in \Gamma(KL)$. Then $\gamma = \alpha+u$ for $x^\alpha\in KL$ and $u\in\R_{\geq0}^m$. The ideal $KL$ is spanned by products $x^{\alpha'}x^{\beta'}$ with $x^{\alpha'}$ a monomial of $K$ and $x^{\beta'}$ a monomial of $L$ as an $\F$-vector space. Since the spanning set consists of monomials and monomials are linearly independent, the monomial $x^\alpha$ is equal to one of them. Let $x^\alpha = x^{\alpha'}x^{\beta'}$ with $x^{\alpha'}\in K, x^{\beta'}\in L$. Therefore, $\gamma = \alpha'+(\beta'+u)\in \Gamma(K)+\Gamma(L)$, and hence $\Gamma(KL)\subset\Gamma(K)+\Gamma(L)$. So (iv) holds.
        \end{proof}

        A monomial order is a well order $\leq$ on the set of monomials in $P$ that satisfies $m m_1\leq m m_2$ if $m_1\leq m_2$ for any monomials $m,m_1,m_2\in P$.
        Fix a monomial ordering $\leq$ on $P$. For a nonzero polynomial $f\in P$, we denote by $\LM(f)$ its greatest term with respect to $\leq$, with its coefficient removed. 

        \begin{lemma}\label{lem: LM}
            Let $f,g\in P$ be nonzero. Then $\LM(fg) = \LM(f)\LM(g)$.
        \end{lemma}
        \begin{proof}
            Every monomial $m_1m_2$ occurring in $fg$ satisfies $m_1m_2\leq \LM(f)m_2\leq \LM(f)\LM(g)$ by the definition of monomial order. If there exists a monomial $n_1n_2$ occurring in $fg$ such that $n_1n_2=\LM(f)\LM(g)$ with $n_1<\LM(f)$, then $n_1n_2<\LM(f)n_2\leq \LM(f)\LM(g) = n_1n_2$, which is a contradiction.
        \end{proof}

        \begin{proposition}\label{prop: inequality}
            Let $K, L$ be two ideals of $P$ that have finite colength. Then
            \begin{equation}\label{inequality}
                \colen(KL)^{\frac{1}{m}} \leq \colen(K)^{\frac{1}{m}} + \colen(L)^{\frac{1}{m}}.
            \end{equation} 
        \end{proposition}
        \begin{proof}
            If one of the ideals $K, L, KL$ is equal to $P$, the inequality is immediate. If one of them is $0$, its colength is infinite, which contradicts the assumption that $K, L, KL$ have finite colength. Assume that $K, L, KL$ are nonzero and proper.

            First assume that $K$ and $L$ are monomial ideals. 
            Let $\Theta = S^{m-1}\cap \R_{\geq 0}^m$ where $S^{m-1}$ is the unit sphere in $\R^m$, and let $\diff \theta$ be the surface measure on $\Theta$. For $\theta\in\Theta$, set
            \begin{equation*}
                \rho_K(\theta) = \min\{t\geq 0 : t\theta\in\Gamma(K)\}.
            \end{equation*}
            By Lemma \ref{lem:inequality} (i), $\Gamma(K)$ is closed. By Lemma \ref{lem: colength}, $K$ contains a power of every variable; then each ray in $\R_{\geq 0}^m$ meets $\Gamma(K)$. Therefore, $\rho_K(\theta)$ is well-defined.       
            The function $\rho_K$ is bounded by $\sqrt{m}\max\{ a_i:1\leq i\leq m\}$, where $x_i^{a_i}\in K$ are as in Lemma \ref{lem: colength}, and it is measurable since $\{\theta:\rho_K(\theta)\leq t\} = \{\theta: t\theta\in\Gamma(K)\}$ is closed for $t\geq 0$. In polar coordinates, by Lemma \ref{lem:inequality} (ii),
            \begin{equation*}
                \colen(K) = \vol(\R_{\geq 0}^m\setminus\Gamma(K)) = \int_{\Theta}\int_{0}^{\rho_K(\theta)} t^{m-1} \diff t\;\diff\theta = \frac{1}{m}\int_{\Theta}\rho_K(\theta)^m\diff \theta.
            \end{equation*}
            Since $\rho_K(\theta)\theta\in\Gamma(K)$ and $\rho_L(\theta)\theta\in\Gamma(L)$, by Lemma \ref{lem:inequality} (iv), their sum is in $\Gamma(KL)$. By Lemma \ref{lem:inequality} (iii), $KL$ has finite colength. Therefore, $\rho_{KL}(\theta) \leq \rho_K(\theta)+\rho_L(\theta)$ for $\theta\in\Theta$. Apply Minkowski's inequality to $L^m(\Theta,\frac{\diff\theta}{m})$, the inequality \eqref{inequality} holds.

            Now let $K$ and $L$ be arbitrary ideals. Let $\ini(K)$ be the initial ideal of $K$, the monomial ideal generated by $\LM(f)$ for any $f\in K$. 
            
            The monomials outside $\ini(K)$ form a basis of $P/K$ under the canonical projection. 
            For independence, let $r$ be a combination of such monomials with $r\in K$. If $r\neq 0$, then $\LM(r)$ is one of them, hence outside $\ini(K)$. But $r\in K$ implies $\LM(r)\in\ini(K)$, which is a contradiction. Therefore, $r= 0$, and the monomials outside $\ini(K)$ are linearly independent.
            For spanning, let $f\in P$ with $\LM(f)\in\ini(K)$. Then $\LM(f) = x^\alpha\LM(h)$ for some nonzero $h\in K$, by Lemma \ref{lem: LM}, $\LM(f) = \LM(x^\alpha h)$. Let $g = x^\alpha h\in K$, the element $f_1 = f-\LC(f)g$, where $\LC(f)$ is the coefficient of $\LM(f)$ in $f$, satisfies $f_1 = 0$ or $\LM(f_1)<\LM(f)$. Iterating, the leading monomial strictly decreases at each step. This terminates since only finitely many monomials are smaller than a given one. The final remainder $r$ is an $\F$-linear combination of monomials outside $\ini(K)$, and $f-r\in K$. Therefore, $\overline{f} = \overline{r}+K$ in $P/K$. Therefore, $\colen(K)=\colen(\ini(K))$.

            By Lemma \ref{lem: LM}, $\ini(K)\ini(L)\subset \ini(KL)$. Therefore,
            \begin{equation*}
                \begin{aligned}
                    \colen(KL)^{\frac{1}{m}} &= \colen(\ini(KL))^{\frac{1}{m}} \leq \colen(\ini(K)\ini(L))^{\frac{1}{m}}\\ 
                    &\leq \colen(\ini(K))^{\frac{1}{m}}+\colen(\ini(L))^{\frac{1}{m}}\\
                    &= \colen(K)^{\frac{1}{m}}+\colen(L)^{\frac{1}{m}}.
                \end{aligned}
            \end{equation*}
        \end{proof}

    \subsection{Maximal-rank argument ingredient}

        \begin{lemma}[Flanders–Atkinson lemma, See {\cite[Lemma 2.2]{MR3255525}}]\label{lem: IVsubsetRv}
            Let $\F$ be an infinite field, and let $V$ be a finite-dimensional vector space over $\F$. For a commutative unital algebra $R\subset\End_\F (V)$, choose $v\in V$ such that $r = \dim_\F(Rv)$ is maximal, and let $I = \Ann_R(v)$. Then 
            \begin{enumerate}[(i)]
                \item $IV\subset Rv$,
                \item $I(Rv) = 0$,
                \item $I^2 = 0$.
            \end{enumerate}
        \end{lemma}
        \begin{proof}
            Let $b_1v,\cdots, b_rv$ be a basis of $Rv$, where $b_1,\cdots,b_r\in R$. Fix $a\in I,w\in V$, and for any $t\in \F$, the vectors
            \begin{equation*}
                b_1(v+tw),\cdots, b_r(v+tw), a(v+tw)
            \end{equation*}
            are contained in $R(v+tw)$ of dimension at most $r$ by the choice of $v$. Therefore, their exterior product vanishes. Since $av = 0$, this gives the polynomial identity
            \begin{equation*}
                (b_1v+tb_1w)\wedge \cdots\wedge(b_rv+tb_rw)\wedge taw = 0.
            \end{equation*}
            The field is infinite, so the coefficient of $t$ is zero, i.e.,
            \begin{equation*}
                b_1v\wedge \cdots\wedge b_rv\wedge aw = 0.
            \end{equation*}
            Since $b_1v,\cdots,b_rv$ are independent, $b_1v\wedge\cdots\wedge b_rv \neq 0$, and hence
            \begin{equation*}
                aw\in\Span\{b_1v,\cdots,b_rv\} = Rv.
            \end{equation*}
            Since the choice of $a,w$ is arbitrary, $IV\subset Rv$, so (i) holds. Condition (ii) holds since commutativity gives $a(bv) = b(av) = 0$ for $b\in R$. By (i), we have $I^2 V = 0$. Since $R$ is an algebra of $\F$-linear maps, $I^2 = 0$, and hence (iii) holds.
        \end{proof}

    \subsection{Proof of Theorem \ref{thm: linear} (i)}
        
        \begin{proof}[Proof of Theorem \ref{thm: linear}(i)]
            Let $V=\F^n$, and let $R = \F[A_1,\cdots,A_m]$. It is easy to see that
            \begin{equation*}
                \dim_\F \F[A_1,\cdots,A_m] = \dim_{\overline{\F}}\overline{\F}[A_1,\cdots,A_m].
            \end{equation*}
            Then we can replace $\F$ by its algebraic closure. Without loss of generality, we assume $\F$ is algebraically closed. So $\F$ is infinite, and hence Lemma \ref{lem: IVsubsetRv} is available.
            
            Fix $m\in\N_+$. We prove this theorem by induction on $n=\dim_\F V$. For $n=1$, the matrices $A_1,\cdots, A_m$ can be regarded as the elements in $\F$. Then $\F[A_1,\cdots, A_m] = \F$ and 
            \begin{equation*}
                \dim_\F \F[A_1,\cdots,A_m] = 1 = (2^m-1)(1-1)+1.
            \end{equation*}
            Therefore, this theorem holds for $n=1$. Assuming that this theorem holds for the positive integers that are strictly less than $n$, we show that this theorem holds for $n$. Choose the vector $v\in V$ such that $r = \dim_\F(Rv)$ is maximal, and let $W = Rv$. Since $R$ has an identity, $r\geq 1$. 

            If $W=V$, then the $\F$-linear map
            \begin{equation*}
                \begin{aligned}
                    \varphi: R &\to V\\
                             a &\mapsto av.
                \end{aligned}            
            \end{equation*}
            is surjective. If $b\in\ker(\varphi)$, then $bv = 0$. Commutativity gives $b(cv) = c(bv) = 0$ for any $c\in R$, i.e., $b V = 0$. Since $R$ is an algebra of $\F$-linear maps, $b=0$. Then $\ker(\varphi) = \{0\}$. It follows that $\varphi$ is also injective. Therefore, $\varphi$ is an isomorphism between $R$ and $V$. So we have
            \begin{equation*}
                \dim_\F R = \dim_\F V = n \leq (2^m-1)(n-1)+1,
            \end{equation*}
            the last inequality holds since $2^m -1\geq 1$ for $m\in\N_+$.

            Assume $\dim_\F V/W>0$. Regard $V$ as a $P=\F[x_1,\cdots,x_m]$-module via the action
            \begin{equation*}
                f(x_1,\cdots,x_m)\cdot v = f(A_1,\cdots,A_m)v.
            \end{equation*}
            Then 
            \begin{equation}\label{eq: R=P/Ann(V)}
                R\cong P/\Ann_P(V).    
            \end{equation}
            Commutativity and $W= Pv$ give $\Ann_P(W) = \Ann_P(v)$, so $f\mapsto fv$ induces $P/\Ann_P(W)\cong W$, and hence 
            \begin{equation}\label{eq: codim(Ann_P(W))=dim(W)=r}
                \colen(\Ann_P(W)) = \dim_\F(W) = r.
            \end{equation} 
            The image of $\Ann_P(W) = \Ann_P(v)$ in $R = P/\Ann_P(V)$ is $\Ann_R(v)$, by Lemma \ref{lem: IVsubsetRv}(i), $\Ann_P(W)V\subset W$, and hence $\Ann_P(W)\subset\Ann_P(V/W)$. It follows that
            \begin{equation}\label{eq: codim(Ann(V/W)leq r)}
                \colen(\Ann_P(V/W))\leq\colen(\Ann_P(W)) = r.
            \end{equation}
            Let $R_{V/W}$ be the image of $P$ in $\End_\F(V/W)$, i.e., the algebra induced by $R$ on $V/W$, generated by the identity and by the quotient map induced by $A_1,\cdots, A_m$.
            The action map $P\to \End_\F(V/W)$ has kernel $\Ann_P(V/W)$ and image $R_{V/W}$, by the rank-nullity theorem and the inductive hypothesis, 
            \begin{equation}\label{eq: thm inductive step}
                \colen(\Ann_P(V/W)) = \dim_\F(R_{V/W}) \leq (2^m-1)(\dim_\F(V/W)-1)+1.
            \end{equation}
            Since $\Ann_P(V/W)$ maps $V$ into $W$ and $\Ann_P(W)$ annihilates $W$, the product $\Ann_P(W)\Ann_P(V/W)$ annihilates $V$, and hence 
            \begin{equation}\label{eq: Ann(W)Ann(V/W)subset Ann(V)}
                \Ann_P(W)\Ann_P(V/W)\subset \Ann_P(V).
            \end{equation}
            Therefore, let $s= \colen(\Ann_P(V/W))$, by Formulas \eqref{eq: R=P/Ann(V)}, \eqref{eq: codim(Ann_P(W))=dim(W)=r}, \eqref{eq: codim(Ann(V/W)leq r)}, \eqref{eq: thm inductive step}, \eqref{eq: Ann(W)Ann(V/W)subset Ann(V)}, and Proposition \ref{prop: inequality}, 
            \begin{equation*}
                \begin{aligned}
                    \dim_\F(R) &= \colen(\Ann_P(V))\\
                    &\leq \colen(\Ann_P(W)\Ann_P(V/W))\\
                    &\leq \left(\colen(\Ann_P(W))^{\frac{1}{m}} + \colen(\Ann_P(V/W))^{\frac{1}{m}}\right)^m\\
                    &= \left(r^{\frac{1}{m}}+s^{\frac{1}{m}}\right)^m\\
                    &= r+s + \sum_{j=1}^{m-1}\binom{m}{j} r^{\frac{m-j}{m}}s^{\frac{j}{m}}\\
                    &\leq r+ (2^m-1)(\dim_\F(V/W)-1)+1 + (2^m-2)r\\
                    &=(2^{m}-1)(r+\dim_\F(V/W)-1)+1\\
                    &=(2^{m}-1)(\dim_\F(W)+\dim_\F(V/W)-1)+1\\
                    &=(2^{m}-1)(n-1)+1,
                \end{aligned}
            \end{equation*}
            where the sum is empty when $m=1$. This completes the induction.
        \end{proof}

    \section{Proof of Theorem \ref{thm: linear}(ii)}

    \subsection{Preliminaries}

        Fix integers $p,q\geq2$ and set $m=p+q$. Let
        \begin{equation*}
            \B = \bigoplus_{r\geq 0}\B_r, \quad\B_r = \F[x_0,\cdots,x_{p-1},z_0,\cdots,z_{q-1}]_{(r,r)},
        \end{equation*}
        where the subscript $(r,r)$ denotes the bi-degree, with degree $r$ with respect to the $x$-variables and degree $r$ with respect to the $z$-variables.

        For $d\geq 0$, let $V_d$ be an $\F$-vector space with the basis
        \begin{equation*}
            v_{\alpha,\beta} \quad (\alpha\in\N^p, \beta\in\N^q, |\alpha|+|\beta| = d).
        \end{equation*}
        Define
        \begin{equation*}
            C_{ij} v_{\alpha,\beta} = \begin{cases}
                                        v_{\alpha+\e_i,\beta-\f_j}, &\beta_j>0,\\
                                        0, &\beta_j=0,
                                    \end{cases}
            \quad
            0\leq i<p,\;0\leq j<q.
        \end{equation*}
        where $\e_i,\f_j$ are the respective coordinate vectors. 
        For $|\mu|=|\nu|$, define an $\F$-algebra homomorphism $\rho_d:\mathcal{B} \to \End_\F(V_d)$ given by 
        \begin{equation*}
            \rho_d(x^\mu z^\nu)v_{\alpha,\beta} = \begin{cases}
                                                        v_{\alpha+\mu,\beta-\nu},&\text{$\beta\geq\nu$ coordinatewise,}\\
                                                        0,&\text{otherwise.}
                                                \end{cases}
        \end{equation*}

        \begin{lemma}\label{lem: rho}
            \begin{enumerate}[(i)]
                \item $\dim_\F V_d = \binom{d+m-1}{m-1}$.
                \item $\{C_{ij} : 0\leq i<p, 0\leq j<q\}$ commute pairwise.
                \item For each $d\geq 0$,
                \begin{equation*}
                    \ker(\rho_d) = \B_{>d} := \bigoplus_{r>d}\B_r.
                \end{equation*}
                \item for any subalgebra $\A = \bigoplus\limits_{r\geq0}\A_r$ of $B$, the map $\rho_d$ induces a faithful action of $\A/\A_{>d}$ on $V_d$.
            \end{enumerate}
        \end{lemma}
        \begin{proof}
            (i) Counting the number of bases of $v_{\alpha,\beta}$ is equivalent to counting the number of tuples $(\alpha_0,\cdots,\alpha_{p-1},\beta_1,\cdots,\beta_{q-1})$ satisfying
            \begin{equation*}
                \alpha_0 + \cdots + \alpha_{p-1} + \beta_0 + \cdots +\beta_{q-1} = d.
            \end{equation*}
            The latter can be counted by the stars-and-bars method; there are $\binom{d+m-1}{m-1}$ in total. Therefore,
            \begin{equation*}
                \dim_\F V_d = \binom{d+m-1}{m-1}.
            \end{equation*}

            (ii) It is easy to see that $C_{ij}C_{kl} = C_{kl}C_{ij}$.

            (iii) Let $x^\mu z^\nu \in \B_r$ with $r>d$, and let $v_{\alpha,\beta}$ be a basis element of $V_d$. Then $|\alpha| + |\beta| = d$. Since 
            \begin{equation*}
                |\nu|=r>d\geq |\beta|,
            \end{equation*}
            there exists some $0\leq i_0\leq q-1$ such that $\beta_{i_0}<\nu_{i_0}$, and hence $x^\mu z^\nu$ annihilates $v_{\alpha,\beta}$. Therefore, $\B_r$ kills $V_d$ for $r>d$, i.e., $\B_{>d}\subset\ker(\rho_d)$.

            For the reverse inclusion, let $P = \sum\limits_{i=0}^s P_i = \sum\limits_{i=0}^{d}P_i+\sum\limits_{j=d+1}^s P_i \in \ker(\rho_d)$ where $P_i\in\mathcal{B}_i$. Then the part $\sum\limits_{j=d+1}^s P_j$ kills $V_d$, which implies $\sum\limits_{j=d+1}^s P_j\in\ker(\rho_d)$.
            For the remainder part $P' = \sum\limits_{i=0}^{d} P_i$, assume, by contradiction, that $P'\neq 0$, i.e., $P'=\sum\limits_{i=0}^{r} P_i$, where $r\leq d, P_r\neq 0$. Let
            \begin{equation*}
                P_r'(x,z) = \sum_{|\nu|=r}P_\nu(x)z^\nu
            \end{equation*}
            and choose $\nu_0$ with $P_{\nu_0}\neq 0$. For any $\gamma\in\N^p$ with $|\gamma|=d-r$, and the basis element $v_{\gamma,\nu}$. Then 
            \begin{equation*}
                \rho_d(P_r') v_{\gamma,\nu_0} = \sum_{|\mu|=r}c_\mu v_{\gamma+\mu,0}\neq 0
            \end{equation*}
            since $|\nu| = |\nu_0| = r$ and $\nu\leq\nu_0$ implies $\nu=\nu_0$. For $i<r$,
            \begin{equation*}
                \rho_d(P'_{i})v_{\gamma,\nu_0} = \sum_{|\nu|=i} c_\nu v_{\gamma+\mu,\nu_0-\nu}
            \end{equation*}
            where $\nu_0-\nu\neq 0$. The terms of index $<r$ contain only basis elements whose $\beta$-part is nonzero. However, the terms of index $r$ contain basis elements with $\beta$-part $0$. Then $\rho_d(P)v_{\gamma,\nu_0}=\rho_d(P')v_{\gamma,\nu_0}\neq 0$, which contradicts the fact that  $P\in\ker(\rho_d)$. It follows that $P'=0$, i.e., $P = \sum\limits_{j=d+1}^s P_j\in \B_{>d}$, and hence $\ker(\rho_d)\subset\B_{>d}$.

            (iv) By (iii), $\A\cap \ker(\rho_d)=\A_{>d}$, so $\ker(\rho_d|_{\A}) = \A_{>d}$, and hence $\rho_d$ induces a faithful action of $\A/\A_{>d}$ on $V_d$. 
        \end{proof}

    \subsection{A Subalgebra of $\mathcal{B}$}

        \begin{lemma}\label{lem: subalgebra}
            Let $m=p+q, s=p+q-2$, and let
            \begin{equation*}
                h_l = \sum_{\substack{0\leq i <p\\0\leq j<q\\i+j=l}}x_iz_j\quad(0\leq l\leq s), \qquad u=x_{p-1}z_{q-2}.
            \end{equation*}
            and 
            \begin{equation*}
                \A = \F[h_0,\cdots,h_s,u].
            \end{equation*}
            Then
            \begin{enumerate}[(i)]
                \item the elements $h_0,\cdots,h_s$ are algebraically independent.
                \item there is a monic polynomial $f\in\F[H_0,\cdots,H_s,U]$ of degree $e = \binom{m-2}{p-1}$ in $U$ and homogeneous of total degree $e$ such that
                \begin{equation*}
                    \A\cong \F[H_0,\cdots,H_s,U]/(f).
                \end{equation*}
                under $H_i\mapsto h_i, U\mapsto u$. Consequently,
                \begin{equation*}
                    \A = \bigoplus_{i=0}^{e-1}\F[h_0,\cdots,h_s]u^i.    
                \end{equation*}
                \item Let $R_d = \rho_d(\A)=\F[\rho_d(h_0),\cdots,\rho_d(h_s),\rho_d(u)]$. Then
                \begin{equation*}
                    R_d\cong\A/\A_{>d},
                \end{equation*}
                and 
                \begin{equation*}
                    \dim_\F R_d = \sum_{j=0}^{e-1}\binom{d-j+m-1}{m-1} = \binom{d+m}{m} - \binom{d+m-e}{m}.
                \end{equation*}
            \end{enumerate}
        \end{lemma}
        \begin{proof}
            (i) Let $a,\gamma_1,\cdots,\gamma_s$ be algebraically independent indeterminates over $\F$ with $a=z_{q-1}x_{p-1}$, and put
            \begin{equation*}
                f(t)=a\prod_{i=1}^{p-1}(t-\gamma_i)=\sum_{i=0}^{p-1}z_{q-1}x_it^i,\qquad
                g(t)=\prod_{j=p}^{s}(t-\gamma_j)=t^{q-1}+\sum_{j=0}^{q-2}z_jt^j ,    
            \end{equation*}
            Since
            \begin{equation*}
                f(t)g(t) = a\prod\limits_{i=1}^{s}(t-\gamma_i),
            \end{equation*}
            where $\deg f(t) = p-1$ and $\deg g(t)=q-1$, comparing the coefficient of $t^l$ gives
            \begin{equation*}
                \sum_{i+j=l}x_iz_j=(-1)^{s-l}a\varepsilon_{s-l}(\gamma_1,\cdots,\gamma_s)    
            \end{equation*}
            for $0\leq l\leq s$, where $\varepsilon_{s-l}$ is the $(s-l)$-th elementary symmetric polynomial. For simplicity, denote $\varepsilon_i(\gamma_1,\cdots,\gamma_s)$ as $\varepsilon_i$. By Artin's Theorem (see \cite[Corollary 3.11]{MR4506896}), 
            the field extension 
            \begin{equation*}
                \F(\gamma_1,\cdots,\gamma_s)/\F(\gamma_1,\cdots,\gamma_s)^{S_n}
            \end{equation*}
            is finite, so they have the same transcendence degree. Since the transcendence degree of $\F(\gamma_1,\cdots,\gamma_s)$ is $s$,  by the fundamental theorem of symmetric polynomials, the transcendence degree of
            \begin{equation*}
                \F(\gamma_1,\cdots,\gamma_s)^{S_n} \cong \F(\varepsilon_1,\cdots,\varepsilon_s)
            \end{equation*}
            is also $s$. Therefore, $\varepsilon_1,\cdots,\varepsilon_s$ are algebraically independent. If $a,\varepsilon_1\,\cdots,\varepsilon_s$ are algebraically dependent over $\F$, there exists a nonzero polynomial $P\in \F[y_0,\cdots,y_s]$ such that
            \begin{equation*}
                \begin{aligned}
                    P(a,\varepsilon_1,\cdots,\varepsilon_s) = \sum_{i=0}^kQ_i(\varepsilon_1,\cdots,\varepsilon_s)a^i = 0
                \end{aligned}    
            \end{equation*}
            where $Q_i\in\F[y_1,\cdots,y_s]$.
            Since $a$ is an indeterminate over $\F$, this equality implies 
            \begin{equation*}
                Q_1(\varepsilon_1,\cdots,\varepsilon_s) = \cdots = Q_s(\varepsilon_1,\cdots,\varepsilon_s)=0,
            \end{equation*}
            Since $\varepsilon_1,\cdots,\varepsilon_s$ are algebraically independent, these polynomials $Q_1,\cdots,Q_s$ are zero polynomial in $\F[y_1,\cdots,y_s]$. It follows that $P = 0$ in $\F[y_0,\cdots,y_s]$, which contradicts that $P$ is nonzero. Therefore, $a,\varepsilon_1,\cdots,\varepsilon_s$ are algebraically independent. By the definition of $h_0,\cdots,h_s$, we have
            \begin{equation*}\label{eq: hi-s}
                h_s = a, \quad h_{s-j} = (-1)^j a \varepsilon_j(\gamma_1,\cdots,\gamma_s).
            \end{equation*}
            Therefore, $h_0,\cdots,h_s$ are algebraically independent.

            (ii) Consider the polynomial
            \begin{equation}\label{eq: defF*}
                F_*(U) = \prod_{\substack{I\subset\{1,\cdots,s\}\\|I|=q-1}}\left(U+a\sum_{i\in I}\gamma_i\right) = U^e + \sum_{i=1}^{e-1} a^iP_i(\gamma_1,\cdots,\gamma_s) U^{e-i},
            \end{equation}
            where
            \begin{equation*}
                P_i(\gamma_1,\cdots,\gamma_s) = \sum_{I_1\neq \cdots\neq I_i\subset \{1,\cdots,s\}} \prod_{t=1}^i\left(\sum_{r\in I_t} \gamma_r\right).
            \end{equation*}
            Since $P_i$ is symmetric and homogeneous of degree $i$, by the fundamental theorem of symmetric polynomials, $P_j(\gamma_1,\cdots,\gamma_s)$ is a linear combination of some $\prod\limits_{j=1}^s\varepsilon_j^{b_j}$ for some $b_j\geq0$ satisfying $i=\sum\limits_{j=1}^s jb_j$. By Formula \eqref{eq: hi-s}, we have
            \begin{equation*}
                \begin{aligned}
                    a^i\prod_{j=1}^s\varepsilon_j^{b_j} &= a^{i-\sum\limits_{j=1}^sb_j}\prod_{j=1}^s\left(a\varepsilon_j\right)^{b_j}\\
                    &= a^{i-\sum\limits_{j=1}^sb_j}\prod_{j=1}^s\left((-1)^jh_{s-j}\right)^{b_j}\\
                    &= a^{i-\sum\limits_{j=1}^s b_j}(-1)^{\sum\limits_{j=1}^sjb_j}\prod_{j=1}^s h_{s-j}^{b_j}\\
                    &= (-1)^i h_s^{i-\sum\limits_{j=1}^sb_j}\prod_{j=1}^s h_{s-j}^{b_j}.
                \end{aligned}
            \end{equation*}
            The exponent of $h_s$ is nonnegative since $i=\sum\limits_{j=1}^sjb_j\geq\sum\limits_{j=1}^sb_j$. Therefore, $F_*(U)$ is the image under the map $H_i\mapsto h_i$ of a polynomial
            \begin{equation*}
                F(H_0,\cdots,H_s,U) = U^e + \sum_{i=1}^{e-1}\widetilde{P}_i(H_0,\cdots,H_{s})U^{e-i}\in\F[H_0,\cdots,H_s,U].
            \end{equation*}
            with $\widetilde{P}_i$ homogeneous of degree $i$. Consider the factor $U+a\sum\limits_{i\in\{p,\cdots,s\}}\gamma_i = U-u$ in Formula \eqref{eq: defF*}, it follows that $F_*(u) = 0$. Therefore, $F(h_0,\cdots,h_s,u) = 0$.

            Regard $\F(h_0,\dots,h_s)$ as a subfield of $\F(a,\gamma_1,\cdots,\gamma_s)$. By Formula \eqref{eq: hi-s}, the field $\F(h_0,\cdots,h_s)$ is generated over $\F$ by $a$ and the $a\varepsilon_i$, and hence also by $a$ and the $\varepsilon_i$ since $\varepsilon_i=(-1)^i\frac{h_{s-i}}{h_s}$. Every permutation of $\gamma_1,\cdots, \gamma_s$ fixes $a$ and each $\varepsilon_i$, so it fixes $K$, and it sends $u$ through the $e$ distinct elements
            \begin{equation*}
                -a\sum_{i\in I}\gamma_i
            \end{equation*}
            for some $|I| = q-1$. Therefore, if $G(U)\in \F(h_0,\cdots,h_s)[U]$ is nonzero and $G(u)=0$, applying these permutations shows that $G$ has at least $e$ distinct roots, so $\deg_U G \geq e$. Since $F_*$ has degree $e$ and is monic, it is the minimal polynomial of $u$ over $\F(h_0,\cdots,h_s)$.

            Let $G(H_0,\cdots,H_s,U)\in\F[H_0,\cdots,H_s,U]$ satisfy $G(h_0,\cdots,h_s,u)=0$. Dividing by the monic polynomial $f$ in the variable $U$, we have
            \begin{equation*}
                G=QF+R,    
            \end{equation*}
            where $Q,R\in\F[H_0,\cdots,H_s][U]$ and $\deg_U R<e$. Then $R(h_0,\cdots,h_s,u)=0$. 
            Assume that $R=\sum\limits_{i=0}^kQ_i(H_0,\cdots,H_s)U^i$ with $Q_i\in\F[H_0,\cdots,H_s]$. If $R\neq 0$, then there exists some $Q_{i_0}\neq 0$. Substitute $H_j\mapsto h_j$, then $R(h_0,\cdots,h_s, U)\in \F(h_0,\cdots,h_s)[U]$ of degree $<e$. By the minimality of the degree of $F$, we have each $Q_i(h_0,\cdots,h_s)=0$. It follows that $Q_{i_0}(h_0,\cdots,h_s)=0$, which contradicts the fact $h_0,\cdots,h_s$ are algebraically independent. Therefore, $G = QF$ for some $Q\in\F[H_0,\cdots,H_s,U]$.  

            (iii) By Lemma \ref{lem: rho}, we have $\ker(\rho_d|_{\A}) = \A_{>d}$. Therefore,
            \begin{equation}\label{eq: Rd}
                R_d\cong \A/\ker(\rho_d|_\A) =  \A/\A_{>d},
            \end{equation}
            and Lemma \ref{lem: rho} show that $R_d$ acts faithfully on $V_d$.

            Let $\F[h_0,\cdots,h_s]_{\leq d-j}$ be the $\F$-subspace of polynomials of degree $\leq d-j$. 
            A monomial of total degree at most $d-j$ in $s+1=m-1$ variables is determined by its exponent vector $\alpha=(\alpha_0,\cdots,\alpha_s)$ such that
            \begin{equation}\label{eq: dimsumleq}
                \alpha_0+\cdots+\alpha_s\leq d-j.
            \end{equation}
            Add one variable $\alpha_{s+1}$ such that
            \begin{equation}\label{eq: dimsumeq}
                \alpha_0+\cdots+\alpha_s+\alpha_{s+1} = d-j,
            \end{equation}
            then counting the number of tuples $(\alpha_0,\cdots,\alpha_{s+1})$ satisfying Formula \eqref{eq: dimsumleq} is equivalent to counting the number of tuples $(\alpha_0,\cdots,\alpha_{s+1})$ satisfying Formula \eqref{eq: dimsumeq}. The latter can be counted by the stars-and-bars method; there are $\tbinom{d-j+m-1}{m-1}$ in total. Therefore,
            \begin{equation}\label{eq: dimF[h]<=d-j}
                \dim_\F\F[h_0,\cdots,h_s]_{\leq d-j} = \binom{d-j+m-1}{m-1}.
            \end{equation}
            By (ii) and $u\in \B_1$
            \begin{equation*}
                \A/\A_{>d} = \bigoplus_{j=0}^{e-1} \F[h_0,\cdots,h_s]_{\leq d-j} u^j,
            \end{equation*}
            by Formulas \eqref{eq: Rd}, \eqref{eq: dimF[h]<=d-j}, we have
            \begin{equation}\label{eq: dimfirst}
                \begin{aligned}
                    \dim_\F R_d &=  \dim_\F \A/\A_{>d} = \sum_{j=0}^{e-1}\dim_\F \F[h_0,\cdots,h_s]_{\leq d-j}\\ &= \sum_{j=0}^{e-1}\binom{d-j+m-1}{m-1}.
                \end{aligned}
            \end{equation}

            Consider the polynomial
            \begin{equation*}
                \begin{aligned}
                    \sum_{j=0}^{e-1}(1+t)^{d-j+m-1}  &=  (1+t)^{d+m-e}((1+t)^0 + \cdots + (1+t)^{e-1})\\
                    &= (1+t)^{d+m-e}\frac{(1+t)^e-1}{t}\\
                    &= \frac{1}{t}(1+t)^{d+m} - \frac{1}{t}(1+t)^{d+m-e}.
                \end{aligned}
            \end{equation*}
            By Formula \eqref{eq: dimfirst}, $\dim_\F R_d$ is the coefficient of the term $t^{m-1}$ on the left-hand side, and on the right-hand side, $t^{m-1}$ has coefficient $\binom{d+m}{m}-\binom{d+m-e}{m}$.
            Therefore,
            \begin{equation*}
                \dim_\F R_d = \binom{d+m}{m}-\binom{d+m-e}{m}.
            \end{equation*}
        \end{proof}

    \subsection{Proof of Theorem \ref{thm: linear}(ii)}

        \begin{proof}[Proof of Theorem \ref{thm: linear}(ii)]
            Fix $p,q$. By Lemma \ref{lem: subalgebra}, there exists an algebra $R_d=\F[\rho_d(h_0),\cdots,\rho_d(h_s),\rho_d(u)]$ such that
            \begin{equation}\label{eq: 2final1}
                \dim_\F R_d = \binom{d+m}{m}-\binom{d+m-e}{m}, 
            \end{equation}
            where $e = \binom{m-2}{p-1}$.
            Then, by Lemma \ref{lem: rho} and Formula \eqref{eq: 2final1},
            \begin{equation*}
                c_m(\F) \geq \frac{\dim_\F R_d}{\dim_\F V_d} = \frac{\binom{d+m}{m}-\binom{d+m-e}{m}}{\binom{d+m-1}{m-1}}.
            \end{equation*}
            Letting $d\to\infty$, we have
            \begin{equation*}
                c_m(\F) \geq \lim_{d\to\infty}\frac{\binom{d+m}{m}-\binom{d+m-e}{m}}{\binom{d+m-1}{m-1}}=e.
            \end{equation*}
            Let $p=\lfloor\frac{m-2}{2}\rfloor+1$ and $q=m-p$. Then
            \begin{equation}\label{eq: 2final2}
                c_m(\F) \geq e = \binom{m-2}{p-1}=\binom{m-2} {\lfloor \frac{m-2}{2} \rfloor}.
            \end{equation}

            For each integer $0\leq i\leq m-2$, we have
            \begin{equation*}
                \binom{m-2} {\lfloor \frac{m-2}{2} \rfloor}\geq \binom{m-2}{i}
            \end{equation*}
            Therefore,
            \begin{equation}\label{eq: 2final3}
                \binom{m-2} {\lfloor \frac{m-2}{2} \rfloor} \geq \frac{1}{m-1}\sum_{i=0}^{m-2}\binom{m-2}{i} = \frac{2^{m-2}}{m-1}
            \end{equation}
            By Theorem \ref{thm: linear} (i), Formulas \eqref{eq: 2final2} and \eqref{eq: 2final3}, we have
            \begin{equation*}
                \sqrt[m]{\frac{2^{m-2}}{m-1}} \leq \sqrt[m]{c_m(\F)} \leq r_m \leq \sqrt[m]{\sup_n\frac{(2^m-1)(n-1)+1}{n}}\leq \sqrt[m]{2^m-1}.
            \end{equation*}
            Letting $m\to\infty$, we have
            \begin{equation*}
                \lim_{m\to \infty}\sqrt[m]{\frac{2^{m-2}}{m-1}} = 2 \quad \text{and}\quad \lim_{m\to \infty}\sqrt[m]{2^m-1} = 2
            \end{equation*}
            and hence
            \begin{equation*}
                \lim_{m\to\infty}\sqrt[m]{c_m(\F)} = \lim_{m\to\infty}\sqrt[m]{r_m} = 2.
            \end{equation*}
        \end{proof}

    \section{Proof of Theorem \ref{thm: linear}(iii)}

    \subsection{Preliminaries}

        \begin{lemma}\label{lem: fieldDescending}
            Let $\K$ be a subfield of a field $\F$, and let $A_1,\cdots,A_m\in\M_n(\F)$ commute pairwise. Then there exists a finite extension $\L/\K$ and pairwise commuting matrices $\widetilde{A_1},\cdots,\widetilde{A_m}\in\M_n(\L)$ such that
            \begin{equation*}
                \dim_\L \L[\widetilde{A_1},\cdots,\widetilde{A_m}] \geq \dim_\F\F[A_1,\cdots,A_m].
            \end{equation*}
        \end{lemma}
        \begin{proof}
            Let $M_1,\cdots,M_d$ be an $\F$-basis of $\F[A_1,\cdots,A_m]$. Regarding each $M_i$ as a column vector of its $n^2$ entries, we obtain an $n^2\times d$ matrix over $\F$ of rank $d$, and hence there exists a nonzero $d\times d$ minor $\Delta$. Let
            \begin{equation*}
                B = \K[\text{all entries of $A_1,\cdots,A_m$}, \{d_{ij}\}, \Delta^{-1}]\subset \F.
            \end{equation*} 
            where $\{d_{ij}\}$ is the set of coefficients of one linear combination of $M_1,\cdots,M_d$ under the monomials generated $A_1,\cdots,A_m$.
            Only finitely many entries occur, so $B$ is a finitely generated $\K$-algebra. Let $\m$ be a maximal ideal of $B$, and let $\L= B/\m$. By the weak version of Hilbert’s Nullstellensatz (see \cite[Corollary 7.10]{MR242802}), $\L$ is a finite extension of $\K$.

            Let $\pi:B\to \L$ be the canonical quotient map, and for each $i$, let $\widetilde{A_i}$ be the matrix whose $(s,t)$-entry is the image under $\pi$ of the $(s,t)$-entry of $A_i$. Since $A_i$ and $A_j$ commute, we have $A_iA_j-A_jA_i=0$. Each entry of the right-hand side is a polynomial in the entries of $A_i$ and $A_j$, and $\pi$ is a ring homomorphism. Therefore, every entry of $\widetilde{A_i}\widetilde{A_j}-\widetilde{A_j}\widetilde{A_i}$ is zero, i.e., $\widetilde{A_i}\widetilde{A_j}-\widetilde{A_j}\widetilde{A_i} = 0$. It follows that $\widetilde{A_1},\cdots,\widetilde{A_m}$ commute pairwise. Since $\Delta\cdot \Delta^{-1}=1$ holds in $B$, we have $\pi(\Delta)\cdot\pi(\Delta^{-1}) = \pi(\Delta\cdot\Delta^{-1})= \pi(1_\F) = 1_\L$. Then $ \pi(\Delta)$ is a unit of $\L$. 
            
            For each $i$, let $\widetilde{M_i}$ be the matrix whose $(i,j)$-entry is the image under $\pi$ of the $(i,j)$-entry of $M_i$. We arrange $\widetilde{M_1},\cdots,\widetilde{M_d}$ in the same way as $M_1,\cdots,M_d$. Then the minor $\widetilde{\Delta}$ corresponding to the position $\Delta$ is the image of $\Delta$ under $\pi$. Since $\widetilde{\Delta} = \pi(\Delta)$ is a unit of $\L$, i.e., $\widetilde{\Delta}\neq 0$, and hence $\widetilde{\Delta}$ is a nonzero $d\times d$ minor. It follows that $\widetilde{M_1},\cdots,\widetilde{M_d}$ are linearly independent. Therefore, we have
            \begin{equation*}
                \dim_\L\L[\widetilde{A_1},\cdots,\widetilde{A_m}] \geq d  = \dim_\F \F[A_1,\cdots,A_m].
            \end{equation*}
        \end{proof}

        \begin{lemma}\label{lem: decomposition}
            Let $\L$ be a field, and let $\widetilde{A_1},\cdots,\widetilde{A_m}\in\M_n(\L)$ commute pairwise with characteristic polynomials $\p_1,\cdots,\p_m$ respectively. Consider the splitting field $\E$ over $\L$ of $\p_1,\cdots,\p_m$, and let $V = \E^n$. Then 
            \begin{enumerate}[(i)]
                \item $\E/\L$ is a finite extension.
                \item the algebra $\E[\widetilde{A_1},\cdots,\widetilde{A_m}]$ can be decomposed as a direct sum of local $\E$-algebras, i.e.,
                \begin{equation*}
                    \E[\widetilde{A_1},\cdots,\widetilde{A_m}] = \bigoplus_i R_i.
                \end{equation*}
                \item $V$ can be decomposed into the corresponding simultaneous generalized eigenspaces
                \begin{equation*}
                    V = \bigoplus_{j} V_j.
                \end{equation*}
                \item each $R_j$ acts faithfully on $V_j$.
            \end{enumerate}             
        \end{lemma}
        \begin{proof}
            (i) Let $P(x)=\prod\limits_{i=1}^m\p_i(x)\in\L[x]$. Then $\E$ is the splitting field of $P$ over $\L$ given by adding finitely many roots of finitely many polynomials. Therefore, $\E/\L$ is a finite extension.

            (ii) The $\E$-algebra $\E[\widetilde{A_1},\cdots,\widetilde{A_m}]$ has finite dimension over $\E$, hence an Artinian ring. By the structure theorem for Artinian rings, it can be decomposed as a direct sum of local $\E$-algebras 
            \begin{equation*}
                \E[\widetilde{A_1},\cdots,\widetilde{A_m}] = \bigoplus_{i=1}^sR_i.
            \end{equation*}

            (iii) Let $e_i$ be the idempotent element of $R_i$. Then
            \begin{equation}\label{eq: idempotent}
                \sum_{i=1}^s e_i = 1, \quad e_ie_j = \begin{cases}
                                                        e_i, & i=j,\\
                                                        0, & i\neq j.
                                                    \end{cases}
            \end{equation}
            and $R_i = \E[\widetilde{A_1},\cdots,\widetilde{A_m}]e_i$.
            Let $V_i = e_i V$, by Formula \eqref{eq: idempotent}, we have 
            \begin{equation*}
                V = \bigoplus\limits_{i=1}^s V_i.
            \end{equation*}
            Since $R_i$ is integral over $\E$, $\E\subset \E[\widetilde{A_i}e_j]\subset R_j$, and hence $R_j$ is integral over $\E[\widetilde{A_i}e_j]$. Since $R_j$ is local, by the going-up theorem (see \cite[Theorem 5.11]{MR242802}), $\E[\widetilde{A_i}e_j]$ has a unique maximal ideal. Let $\m_{ij}(x)$ be the minimal polynomial of $\widetilde{A_i}e_j$ on $V_j$. Then
            \begin{equation*}
                \E[\widetilde{A_i}e_j] \cong \E[x]/(\m_{ij}(x)\E[x]).
            \end{equation*}
            Since the left-hand side is local, there exists a unique maximal ideal of the right-hand side. This implies that $\m_{ij}(x) = p(x)^t$ for some irreducible polynomial $p(x)$. Let $\lambda$ be a root of $p(x)$. Then $\lambda$ is a root of $\m_{ij}(x)$, i.e, $\lambda$ is the eigenvalue of $\widetilde{A_i}e_j$ on $V_j$, and hence $\lambda\in\E$. Therefore, $p(x) = (x-\lambda)^k$, and hence $\m_{ij}(x)=(x-\lambda)^{kt}$. It follows that $\widetilde{A_i}e_j$ has only one eigenvalue on $V_j$.

            (iv) Assume that $Ae_i\in R_i$ acts as zero in $V_i$ for some $A\in\E[\widetilde{A_1},\cdots,\widetilde{A_m}]$. 
            For any $j=1,\cdots,s$, 
            $Ae_ie_jv = A(0)=0$ if $i\neq j$, $Ae_ie_jv = Ae_iv = 0$ if $i=j$. So $Ae_i$ acts as zero on $V$. Since $\E[\widetilde{A_1},\cdots,\widetilde{A_m}]$ acts faithfully on $V$, and $Ae_i\in \E[\widetilde{A_1},\cdots,\widetilde{A_m}]$, hence $Ae_i=0$. Therefore, $R_i$ acts faithfully on $V_i$. 
        \end{proof}

    \subsection{Proof of Theorem \ref{thm: linear} (iii)}
       
        \begin{proof}[Proof of Theorem \ref{thm: linear} (iii)]
            Let $\K$ be the prime subfield of $\F$.
            
            For any pairwise commuting matrices $A_1,\cdots,A_m\in\M_n(\K)$, it is easy to see
            \begin{equation*}
                \dim_\K \K[A_1,\cdots,A_m] = \dim_{\F} \F[A_1,\cdots,A_m]\leq c_m(\F).
            \end{equation*}
            Since the choice of $A_1,\cdots, A_m$ is arbitrary, taking the supremum gives $c_m(\K)\leq c_m(\F)$.

            For the reverse inequality, let $A_1,\cdots,A_m\in \M_n(\F)$ commute pairwise. By Lemma \ref{lem: fieldDescending}, there exists a finite extension $\L/\K$ and pairwise commuting matrices $\widetilde{A_1},\cdots,\widetilde{A_m}\in\M_n(\L)$ such that
            \begin{equation}\label{eq: LgeqK}
                \dim_\L \L[\widetilde{A_1},\cdots,\widetilde{A_m}]\geq \dim_\F \F[A_1,\cdots,A_m].
            \end{equation}
            Let $\E$ be the splitting field over $\L$ of the characteristic polynomials of $\widetilde{A_1},\cdots,\widetilde{A_m}$. Since $\L/\K$ is finite, so is $\E/\K$.
            By Lemma \ref{lem: decomposition}, there exist decompositions
            \begin{equation*}
                \E[\widetilde{A_1},\cdots,\widetilde{A_m}] = \bigoplus_i R_i, \quad V=\bigoplus_j V_j
            \end{equation*}
            such that $R_j$ acts faithfully on $V_j$. It follows that
            \begin{equation}\label{eq: dimension}
                \frac{\dim_\E \E[\widetilde{A_1},\cdots,\widetilde{A_m}]}{n} = \sum_j\frac{\dim_\E R_j}{n} = \sum_j \frac{\dim_\E V_j}{n}\frac{\dim_\E R_j}{\dim_\E V_j}.
            \end{equation}
            and $R_j$ can be regarded as the subalgebra of $\End_\E(V_j)$ for each $j$. Since $\frac{\dim_\E V_j}{n}\leq 1$, by Formula \eqref{eq: dimension}, there exists at least one component $W = V_{j_0}$ such that
            \begin{equation}
                \frac{\dim_\E R_{j_0}}{\dim_\E W} \geq \frac{\dim_\E\E[\widetilde{A_{1}},\cdots,\widetilde{A_m}]}{n},
            \end{equation}
            otherwise the right-hand side is strictly less than the left-hand side in Formula \eqref{eq: dimension}.
            Let $\lambda_{i}$ be the eigenvalue of $\widetilde{A_i}e_{j_0}|_W$ in $W$, and let $N_i=\widetilde{A_i}e_{j_0}|_{W}-\lambda_i$. Then $N_1,\cdots, N_m$ are pairwise commuting nilpotent operators on $W$, and $R_{j_0} = \E[N_1,\cdots,N_m]$.

             Since $\E/\L$ and $\L/\K$ are finite extensions, the field extension $\E/\K$ is finite. The field $\K$ is the prime subfield of $\F$, so $\K$ is perfect, and hence the finite extension $\L/\K$ is separable. By the primitive element theorem, there exists an element $\alpha\in\E$ such that $\E = \K(\alpha)$. Regard $W$ as an $\K$-vector space, define
            \begin{equation*}
                B_1 = \alpha I + N_1, \qquad B_i = N_i \quad(2\leq i\leq m).
            \end{equation*}
            These $m$ operators commute pairwise. The part $\alpha I$ is semisimple over $\K$ since its minimal polynomial is separable, and it commutes with $N_1$. By Jordan-Chevalley decomposition, the uniqueness of the decomposition gives that the semisimple part of $B_1$ is $\alpha I$, and hence there exists a polynomial $g(x)\in\K[x]$ such that $\alpha I = g(B_1)$. Therefore, 
            \begin{equation*}
                \alpha I\in\K[B_1]\quad\text{and}\quad\K[B_1,\cdots,B_m] = \E[N_1,\cdots,N_m]\subset\End_\K(W).
            \end{equation*}
            The action is faithful since restriction of scalars does not identify distinct endomorphisms.             
            Therefore, by Formula \eqref{eq: LgeqK}, we have
            \begin{equation*}
                \begin{aligned}
                    c_m(\K) &\geq\frac{\dim_\K\K[B_1,\cdots,B_m]}{\dim_\K W}\\
                    &= \frac{[\E:\K]\dim_\E\E[N_1,\cdots,N_m]}{[\E:\K]\dim_\E W}\\
                    &= \frac{\dim_\E\E[N_1,\cdots,N_m]}{\dim_\E W}\\
                    &=\frac{\dim_\E R_{j_0}}{\dim_\E V_{j_0}}\\
                    &\geq \frac{\dim_\E\E[\widetilde{A_1},\cdots,\widetilde{A_m}]}{n}\\
                    &= \frac{\dim_\L\L[\widetilde{A_1},\cdots,\widetilde{A_m}]}{n}\\
                    &\geq \frac{\dim_\F\F[A_1,\cdots,A_m]}{n}.
                \end{aligned}
            \end{equation*}
            Since the choice of $A_1,\cdots, A_m$ is arbitrary, taking the supremum gives $c_m(\F)\leq c_m(\K)$.
        \end{proof}

    \section*{}

    \subsection*{Acknowledgement.}
         We would like to thank Prof. Kaiming Zhao for his encouragement and helpful suggestions, as well as for informing us of an interesting and beautiful recent work \cite{Zhao2026}. The third author is very grateful to Prof. Genkai Zhang for his kind hospitality and encouragement. We would also like to thank the referee for his/her detailed comments and kindness.

    \subsection*{AI disclosure.}
        During the preparation of this manuscript, the authors used LLMs to assist in drafting and checking parts of the mathematical arguments, organizing references, and improving the clarity and readability of the exposition. All mathematical statements and proofs in the final manuscript were written, reviewed, and verified by the authors, who take full responsibility for the content of this manuscript.

    \bibliographystyle{plain}
    \bibliography{Bibliography}
\end{document}